\documentclass{amsart}
\usepackage[procnames]{listings}
\usepackage{color}
\usepackage{cite}
\usepackage{ amssymb }
\numberwithin{equation}{section}
\usepackage{amsmath}
\usepackage{mathtools}
\usepackage{amssymb}
\title{Diameter of  Ramanujan Graphs and  Random Cayley Graphs}
\author{Naser T Sardari}
\date{March 2017}

	\newtheorem{thm}{Theorem}[section]

	\newtheorem{rem}[thm]{Remark}

	\newtheorem{cor}[thm]{Corollary}
	
	\theoremstyle{defi}

	\theoremstyle{pf}

	\numberwithin{equation}{section}

\begin{document}
\maketitle
\begin{abstract} We study the diameter of LPS Ramanujan graphs $X_{p,q}$. We show that the diameter of the bipartite Ramanujan graphs is greater than $  (4/3)\log_{p}(n) +O(1)$ where $n$ is the number of vertices of  $X_{p,q}$. 
We also construct an infinite  family of  $(p+1)$-regular LPS Ramanujan graphs $X_{p,m}$ such that the diameter of these graphs is greater than or equal to $ \lfloor (4/3)\log_{p}(n) \rfloor$. On the other hand, for any $k$-regular Ramanujan graph we show that the distance of only a tiny fraction of all pairs of vertices is greater than $(1+\epsilon)\log_{k-1}(n)$.  We also have some numerical experiments for LPS Ramanujan graphs and random Cayley graphs which suggest that the diameters are asymptotically $(4/3)\log_{k-1}(n)$ and $\log_{k-1}(n)$, respectively.     \end{abstract}

\tableofcontents

\section{Introduction}
\subsection{Motivation}
 The diameter of any $k$-regular graph with $n$ vertices is bounded from below by  $\log_{k-1}(n)$  and it  could get as large as a scalar multiple of  $n$. It is known that the diameter of  any $k$-regular Ramanujan graph is bounded from above by  $2(1+\epsilon)\log_{k-1}(n)$ \cite{Sarnak2}. Lubotzky, Phillips and Sarnak constructed a family of $(p+1)-$regular Ramanujan graphs $X_{p,q}$ \cite{Sarnak2}, where $p$ and $q$ are prime numbers and $q\equiv 1$ mod 4. $X_{p,q}$ is a $p+1$-regular bipartite or non-bipartite graph depending on $p$ being a non-quadratic or quadratic residue mod $q$, respectively. Their construction can be modified for every integer $q$; see \cite{sarnakbook} or \cite{Lubotzky}.    It was expected  that the diameter of the LPS Ramanujan graphs to be bounded from above by $(1+\epsilon) \log_{k-1}(n)$;  see \cite[Chapter~3]{Sarnak}.  However, we show that the diameter of an infinite family of  $p+1$-regular LPS Ramanujan graphs is greater than or equal to
\begin{equation}
\lfloor (4/3)\log_{p}(n)\rfloor.
\end{equation}


 While there are points $x$ and $y$ whose distance is large in a LPS Ramanujan graph, we prove that the distance of a tiny fraction of vertices in any $k$-regular Ramanujan graph $G$  is less than $(1+\epsilon)\log_{k-1}(n)$.  In other words, the essential diameter is asymptotic to
 
 $$(1+\epsilon)\log_{k-1}(n),$$  
 where the essential diameter of a graph is $d$ if $99\%$ of the distance of pairs of vertices is less than $d$.  In fact, we prove a stronger result, we show that for every vertex $x$ in a $k$-regular Ramanujan graph $G$ the number of points which cannot be visited by exactly $l$ steps, where $l>(1+\epsilon)\log_{k-1}(n)$, is less than $n^{1-\epsilon}$. So the number is exponentially decaying. In particular, it also recovers $2(1+\epsilon) \log_{k-1}(n)$ as an upper bound on the diameter of $k$-regular Ramanujan graph.  Furthermore,  we give some numerical datas for two families of 6-regular graphs. The first family of graphs are the 6-regular LPS Ramanujan graphs and we denote them by $X_{5,q}$. 
 The second family are the 6-regular random Cayley graphs ${\rm PSL}_2(\mathbb{Z}/q\mathbb{Z})$, i.e.  the Cayley graphs that are constructed by 3 random generators of ${\rm PSL}_2(\mathbb{Z}/q\mathbb{Z})$ and their inverses $\{s_1^{\pm},s_2^{\pm},s_3^{\pm}\}$. We denote these graphs by $Z^q$. The numerical experiments suggest that the diameter of the (number theoretic) LPS Ramanujan Graphs is asymptotic to
\begin{equation}
(4/3)\log_{5}(n).
\end{equation}
This is consistent with our conjecture on the optimal strong approximation for quadratic forms in 4 variables \cite{Naser}.
On the other hand, the numerical data suggests that the diameter of the random Cayley graph equals that of a random 6-regular graph \cite{Bolbol}, that is
\begin{equation}
\log_{5}(n).
\end{equation}
The archimedean analog of our question has been discussed in   Sarnak's letter to Scott Aaronson and Andy Pollington; see~\cite{Sarnak3}. In that context, the approximation of points on the  sphere by words in LPS generators is considered. 
This question is related to the theory of quadratic Diophantine equations; see \cite{Naser}. Sarnak defines the notion of the covering exponent and the almost all covering exponent  \cite[Page 3]{Sarnak3} that are the analogue of diameter and the essential diameter in our paper. Sarnak showed that the almost all covering exponent is 1; see \cite[Page 28]{Sarnak3}. Our Theorem~\ref{kk}  is the $p$-adic analogue of Sarnak's theorem. In a recent paper \cite{Peres},  Lubetzky and Peres show the simple random walk exhibits cutoff on Ramanujan Graphs. As a result they give a more detailed version of our Theorem~\ref{kk}. In a similar work, for the family of LPS bipartite Ramanujan graphs, Biggs and Boshier determined the asymptotic behavior of the girth of these graphs; see \cite{Biggs}.  They showed that the girth is asymptotic to 
$$(4/3)\log_{k-1}(n).$$

\subsection{Statement of results}
%
We begin by a brief description of  LPS Ramanujan graphs; see \cite{Sarnak2} for a comprehensive treatment of them. The idea of the construction is coming form number theory, i.e. generalized Ramanujan conjecture. More precisely, we consider the symmetric space ${\rm PGL}_2(\mathbb{Q}_p)/ {\rm PGL}_2(\mathbb{Z}_p)$ which can be identified with a regular $(p+1)$-infinite tree. We note that ${\rm PGL}_2(\mathbb{Z}[1/p])$ acts  from the left on ${\rm PGL}_2(\mathbb{Q}_p)/ {\rm PGL}_2(\mathbb{Z}_p)$. The generalized Ramanujan conjecture, which is a theorem in this case, implies that  the quotient of  ${\rm PGL}_2(\mathbb{Q}_p)/ {\rm PGL}_2(\mathbb{Z}_p)$ by  any congruence subgroup of ${\rm PGL}_2(\mathbb{Z}[1/p])$ which is a $p+1$-regular graph is a Ramanujan graph. By considering an appropriate congruence subgroup of ${\rm PGL}_2(\mathbb{Z}[1/p])$  we can identify the quotient of this symmetric space with a Cayley graph. The Cayley graphs is associated to  ${\rm PSL}_2(\mathbb{Z}/q\mathbb{Z})$ or ${\rm PGL}_2(\mathbb{Z}/q\mathbb{Z})$ 
   depending on $p$ being a quadratic residue or non-quadratic residue mod $q$ where $q$ is a prime and $q\equiv 1$ mod 4. These are LPS Ramanujan graphs that are defined in  \cite{Sarnak2}.

  In the bipartite case where  $p$ is non-quadratic residue mod prime number $p^{\prime}\equiv 1 $ mod 4, we show that the diameter is greater than 
 $$(4/3)\log_{p}(|X_{p,p^{\prime}}|)+O(1).$$
In the non-bipartite case, our theorem is weaker. We take a composite number $m$ such that $p$ is a quadratic residue mod $m$. We show that the diameter of  the LPS Ramanujan graphs $X_{p,m}$ is greater than $$(4/3)\log_{p}(n)-4\log_p (\frac{m}{q})+O(1),$$ where $n=(m^3-m)/2$ is the number of vertices of $X_{p,m}$, $q|m$,  $q$ is a prime power and $q\neq m$. Note that these graphs are congruence covers of the LPS Ramanujan graphs $X_{p,q}$.

In what follows, we give an explicit description of the LPS Ramanujan graphs in terms of the Cayley graphs of ${\rm PSL}_2(\mathbb{Z}/m\mathbb{Z})$. Assume that $q$ is a prime number and $q|m$ where $m$ is an integer and $-1$ is quadratic residue mod $m$. Let $p$ be a prime number such that $p\equiv 1 \text{ mod } 4$ and  $p$ is quadratic residue mod $m$. We denote the  representatives of  square roots of $-1$ and $p$ mod $m$ by $i$ and $\sqrt{p}$, respectively. We are looking at the integral solutions $\alpha=(x_0,x_1,x_2,x_3)$ of the following diophantine equation 
\begin{equation}\label{gen}
x_0^2+x_1^2+x_2^2+x_3^2=p,
\end{equation}
where $x_0 > 0$ and is odd and $x_1,x_2,x_3$ are even numbers. There are exactly $p+1$ integral solutions with such properties . 
%
%
 To each such integral solution $\alpha$, we associates the following matrix $\alpha$ in ${\rm PSL}_2(\mathbb{Z}/m\mathbb{Z})$ :
\begin{equation}
\alpha:=\frac{1}{\sqrt{p}} \begin{bmatrix}
x_0+i x_1 & x_2+i x_3
\\
-x_2+ix_3 & x_0 -ix_1
\end{bmatrix}.
\end{equation}
Note that $\frac{1}{\sqrt{p}}$ is there to make $\text{det}(\alpha)=1$. If $p$ is non-quadratic residue mod $m$ then $\alpha:= \begin{bmatrix}
x_0+i x_1 & x_2+i x_3
\\
-x_2+ix_3 & x_0 -ix_1
\end{bmatrix} \notin {\rm PSL}_2(\mathbb{Z}/m\mathbb{Z}) $ and that's why the Cayley graph in this case is defined over ${\rm PGL}_2(\mathbb{Z}/m\mathbb{Z}) $ and the associated Cayley graph is a bipartite graph.
 This gives us $p+1$ matrices in ${\rm PSL}_2(\mathbb{Z}/m\mathbb{Z})$.  Lubotzky \cite[Theorem 7.4.3]{Lubotzky} showed that they generate ${\rm PSL}_2(\mathbb{Z}/m\mathbb{Z})$ and the associated Cayley graph is a non-bipartite Ramanujan graph. The construction for the bipartite LPS Ramanujan graphs is similar.  The only difference is that $p$ is non-quadratic residue mod $q$.   Furthermore, Lubotzky showed that
%
%

\begin{itemize}
\item diam $X_{p,m} \leq 2 \log_{p} (n)+ 2\log_p 2 +1$.
\item girth $X_{p,m} \geq \frac{2}{3} \log_{p}(n)-2\log_p2$.
\end{itemize}
Our first theorem shows that the distance between the identity matrix $I$ and $W:=\begin{bmatrix} 0 & 1 \\ -1 & 0 \end{bmatrix}$ in the bipartite Ramanujan graph $X_{p,p^{\prime}}=PGL_2(\mathbb{Z}/p^{\prime}\mathbb{Z})$  where $p$ is a non-quadratic residue mod $p^{\prime}$ is bigger than  
$$(4/3)\log_{p}(|X_{p,p^{\prime}}|)-(\frac{2}{3}) \log_{p}2.$$
where $|X_{p,p^{\prime}}|$ is the number of vertices of $X_{p,p^{\prime}}$. For non-bipartite graphs $X_{p,m}$, our theorem shows that either the distance between the identity matrix $I$ and $I^{\prime}:=\begin{bmatrix} 1 & q \\ 0 & 1 \end{bmatrix}$  or between $I$ and $W:=\begin{bmatrix} 0 & 1 \\ -1 & 0 \end{bmatrix}$ in $X_{p,m}$ 
is larger than $$(4/3)\log_{p}(|X_{p,m}|)-4\log_p (\frac{m}{q})-(\frac{2}{3}) \log_{p}2,$$
 As a result,
\begin{equation}
 (4/3)\log_{p}(|X_{p,m}|)-4\log_p (\frac{m}{q})-(\frac{2}{3}) \log_{p}2 \leq \text{diam}(X_{p,m}).
\end{equation} 
\\
%
%

%

%
\begin{thm}\label{ex}\normalfont
 Let $p$, $p^{\prime}$, $q$ be primes and $m$ an integers  defined as above. Let $X_{p,p^{\prime}}$ and $X_{p,m}$ be the associated bipartite and non-bipartite Ramanujan graphs.  Then the diameter of the bipartite LPS Ramanujan graph $X_{p,p^{\prime}}$ is larger than
$$(4/3)\log_{p}(|X_{p,p^{\prime}}|)-(\frac{2}{3}) \log_{p}2.$$ 
In the non-bipartite case, the diameter of the LPS graph $X_{p,m}$ is larger than
 \begin{equation}
 (4/3)\log_{p}(|X_{p,m}|)- 4\log_p (\frac{m}{q})-(\frac{2}{3}) \log_{p}2 \leq \text{diam}(X_{p,m}).
 \end{equation}
\\
\end{thm}
\begin{cor}\normalfont
Let $p$ and $q$ be prime numbers that are congruent to $1$ mod $4$ and $p>1250$. Then the diameter of the LPS Ramanujan graph $X^{p,5q^k}$ for any integer $k$ is greater than or equal to
\begin{equation}
\lfloor (4/3)\log_{p}|X^{p,5q^k}|  \rfloor.
\\
\end{equation} 

\end{cor}
\begin{rem}\normalfont
We conjecture that the diameter of  LPS Ramanujan graph $X_{p,q}$ where $q$ is a prime number is asymptotic to $(4/3)\log_{p}|X_{p,q}|$. I expect  that a variate of our argument gives a sharp lower bound for the diameter of $X_{p,q}$ by choosing  vertices with large distance from the identity  (e.g. $W$ and $I^{\prime}$ in our argument).   We give our numerical results for the distance of $W$ from the identity vertex in Table~\ref{alg}. Our data comes from our algorithm that we developed and implemented  for navigation on LPS Ramanujan graphs \cite{Naser2}. We refer the reader to \cite[Remark 1.10]{Naser2} for further discussion of the distribution of the distance of diagonal elements from the identity where the possible vertices with large distances from the identity matrix are listed. 
\end{rem}

On the other hand, we use the Ramanujan bound on the nontrivial eigenvalues of the adjacency matrix to prove the distance of almost all pairs of vertices is less than $(1+\epsilon)\log_{k}(n)$. The archimedean version of this problem has been discussed in  Sarnak's letter to Scott Aaronson and Andy Pollington \cite[Page~28]{Sarnak3}. More precisely, we prove the following stronger result in Section~\ref{visit}:
\begin{thm}\label{kk}\normalfont
Let $G$ be a $k$-regular Ramanujan graph and fix a vertex $x\in V(G)$. Let $R$ be an integer such that $R>(1+\epsilon)\log_{k-1}(n)$. Define $M(x,R)$ to be the set of all vertices $y\in G$ such that there is no path from $x$ to $y$ with length $R$ (we allow to pass from an edge multiple times but not immediately after). Then,
\begin{equation}
|M(x,R)|\leq  n^{1-\epsilon}(1+R)^2.
\end{equation}

\end{thm}

\subsection{Outline of the paper}
In Section~\ref{mainn}, we prove Theorem~\ref{ex}. The proof  uses some elementary  facts in Number Theory. In Section~{\ref{visit}}, we prove Theorem~\ref{kk}. As a corollary, we prove that the distance of almost all pairs of vertices is less than $$  (1+\epsilon)\log_{k-1}(n).$$
 We use the Chebyshev's inequality by giving an upper bound for  the variance of the distance. We use the Ramanujan bound on the eigenvalues of the adjacency matrix of the graph to give an upper bound on the variance of the distance. Finally in Section~\ref{data}, we compute the diameter of two families of 6-regular graphs. From our numerical experiments, we expect that the diameter of the LPS Ramanujan graphs \cite{Sarnak2}  is asymptotic to
\begin{equation}
(4/3)\log_{p}(n).
\end{equation}
We define a random 6 regular Cayley graph $Z^{q}$ , by considering the Cayley graph of ${\rm PSL}_2(\mathbb{Z}/q\mathbb{Z})$ relative to the generating set $S=\{s_1^{\pm}, s_2^{\pm}, s_3^{\pm} \}$, where $s_1,s_2,s_3$ are random elements of ${\rm PSL}_2(\mathbb{Z}/q\mathbb{Z})$. From the numerical experiments , we  show that in fact the random Cayley graph has a shorter diameter and break the $4/3 \log_5 n$ lower bound for the LPS Ramanujan graphs. For example, We obtained a sample from the random Cayley graph of ${\rm PSL}_2(\mathbb{Z}/229\mathbb{Z})$, such that
\begin{equation}
\text{ diam }(Z^{229}) < 1.23 \log_5 n.
\end{equation}
We expect that the diameter of the random Cayley graph would be as small as possible. More precisely, for $\epsilon >0 $
\begin{equation}
\text{diam}(Z^q) \leq (1+\epsilon) \log_{5}(n)+c_{\epsilon},   \text{ almost surely as } q \to \infty.
\end{equation}

\subsection{Acknowledgments}
I would like to  thank my Ph.D. advisor, Peter Sarnak for suggesting this project to me and also his comments on the earlier versions of this work. I am also very grateful for several insightful and inspiring conversations with him during the course of this work. In addition,   I would like to thank Ori Parzanchevski for sharing his code on computing the largest nontrivial eigenvalue of a Cayley graph with me and also for finding a non-Ramanujan double cover of LPS Ramanujan graphs. Finally, I would like to thank  the careful reading and comments of the anonymous referees. 
\\
\\

\section{Lower bound for the diameter of the Ramanujan graphs}\label{mainn}

In the rest of this section, we give a proof of Theorem \ref{ex}. 
\begin{proof} We begin by proving the first part of the theorem. We show that the distance between the identity matrix $I$ and $W:=\begin{bmatrix} 0 & 1 \\ -1 & 0 \end{bmatrix}$ in the bipartite Ramanujan graph $X_{p,p^{\prime}}=PGL_2(\mathbb{Z}/p^{\prime}\mathbb{Z})$  where $p$ is a non-quadratic residue mod $p^{\prime}$ is bigger than  
$$(4/3)\log_{p}(|X_{p,p^{\prime}}|)-(\frac{2}{3}) \log_{p}2.$$
 By using $|X_{p,p^{\prime}}|=(p^{\prime 3}-p^{\prime})/2$ the above expression simplifies to $4/3 \log (\frac{p^{\prime 3}-p^{\prime}}{4})$ that is smaller than
$$\log_p \frac{p^{\prime 4}}{4}.$$
We proceed by assuming the contradiction that  $dist(I,I^{\prime}) <  \log_p \frac{p^{\prime 4}}{4}$.
There is a  correspondence between non-backtracking path of length $k$ from the identity vertex to another vertex $v_k$ of LPS Ramanujan graph $X_{p,p^{\prime}}$ and the primitive elements of integral quaternion Hamiltonian (the $\gcd$ of the coordinates is one) of square norm $p^k$ up to units of $H(\mathbb{Z})$; see \cite[Section 3]{Sarnak2}. As a result, $dist(I,I^{\prime}) <  \log_p \frac{p^{\prime^4}}{4}$ gives us a solution to the following diophantine equation
\begin{equation}
a^2+b^2+c^2+d^2=p^k,
\end{equation} 
where $k=dist(I,I^{\prime})  $, $b\equiv c \equiv d \equiv 0 \text{ mod } 2p^{\prime}$ and $a \equiv 1 \text{ mod } 2$. At least one of $b, c, d $ is nonzero.  From this we deduce that 
\begin{equation}\label{q2}
a^2 \equiv p^k \text{ mod } p^{\prime 2} \text{ and } 4p^{\prime 2} \leq p^k.
\end{equation}
 Since $p$ is non-quadratic residue mod $p^{\prime}$ the above congruence identity holds only for even  $k$. If $k$ is even and $k=2t$. From \ref{q2} we deduce that
\begin{equation}
a\equiv \pm p^t \text{ mod } p^{\prime 2}.
\end{equation}
If $p^t \geq \frac{p^{\prime 2}}{2} $, 
\begin{equation}
dist(I,I^{\prime})=2t \geq \log_p \frac{p^{\prime 4}}{4},
\end{equation}
a contradiction. Consequently, $p^t < \frac{p^{\prime 2}}{2} $. Since $a\neq \pm p^t$, we deduce that 
\begin{equation}
a=\pm p^t+lp^{\prime 2} \text{  for } l\neq 0.
\end{equation}
Therefore
\begin{equation}
|a|\geq \frac{1}{2}p^{\prime 2}.
\end{equation}
Hence,
\begin{equation}\begin{split}
p^{2t} \geq \frac{ p^{\prime 4}}{4},\text{  and so }
\\
\text{ dist }(I,I^{\prime})=2t \geq \log_p \frac{p^{\prime 4}}{4}.
\end{split}
\end{equation}
a contradiction. Hence, we conclude the first part of our theorem. 
Next we give the proof of the second part of our theorem. Recall that $W:=\begin{bmatrix} 0 & 1 \\ -1 & 0 \end{bmatrix}$  and $I^{\prime}:=\begin{bmatrix} 1 & q \\ 0 & 1 \end{bmatrix}$. By using $|X_{p,m}|=\frac{(m-1)m(m+1)}{2}$ the expression $(4/3) \log_p( n)-4\log_p (\frac{m}{q})-(\frac{2}{3})\log_p2$ simplifies and it is smaller  than 
$$\log_p \frac{q^4}{4}.$$
We  show that 
\begin{equation}
\max(\text{ dist }(I,I^{\prime}),\text{ dist }(I,W)) \geq  \frac{q^4}{4}.
\end{equation}
Assume the contrary that 
\begin{equation}
\max(\text{ dist }(I,I^{\prime}),\text{ dist }(I,W)) <  \frac{q^4}{4}.
\end{equation}
 Similarly by using the correspondence between non-backtracking path of length $k$ and the solutions to the associated  diophantine equation for sums of four squares, $dist(I,I^{\prime}) <  \log_p \frac{q^4}{4}$ gives us a solution to the following diophantine equation
\begin{equation}
a^2+b^2+c^2+d^2=p^k,
\end{equation} 
where $k=dist(I,I^{\prime})$, $b\equiv c \equiv d \equiv 0 \text{ mod } 2q$ and $a \equiv 1 \text{ mod } 2$. At least one of $b, c, d $ is nonzero.  From this we deduce that 
\begin{equation}\label{q2}
a^2 \equiv p^k \text{ mod } q^2 \text{ and } 4q^2 \leq p^k.
\end{equation}
We consider two cases: $k$ even and $k$ odd.

If $k$ is even and $k=2t$. From \ref{q2} we deduce that
\begin{equation}
a\equiv \pm p^t \text{ mod } q^2.
\end{equation}
If $p^t \geq \frac{q^2}{2} $, 
\begin{equation}
dist(I,I^{\prime})=2t \geq \log_p \frac{q^4}{4},
\end{equation}
a contradiction. Consequently, $p^t < \frac{q^2}{2} $. Since $a\neq \pm p^t$, we deduce that 
\begin{equation}
a=\pm p^t+lq^2 \text{  for } l\neq 0.
\end{equation}
Therefore
\begin{equation}
|a|\geq \frac{1}{2}q^2.
\end{equation}
Hence,
\begin{equation}\begin{split}
p^{2t} \geq \frac{ q^4}{4},\text{  and so }
\\
\text{ dist }(I,I^{\prime})=2t \geq \log_p \frac{q^4}{4}.
\end{split}
\end{equation}
a contradiction. Hence $k$ is odd and $k=2t+1$.

We want to use a similar argument to show that $\text{dist } (I,W)= 2t_0+1$ is an odd number. $\text{dist }(I,W) < 4/3 \log_p(n)$ gives us a solution to the following diophantine equation
\begin{equation}
a^2+b^2+c^2+d^2=p^k.
\end{equation} 
Where $b\equiv a \equiv d \equiv 0 \text{ mod } q$ and $c \equiv 0 \text{ mod } 2$. Since $a$ is odd, then $q \leq |a|$.  We deduce that 
\begin{equation}\label{qee}
c^2 \equiv p^k \text{ mod } q^2 \text{ and } q^2 \leq p^k.
\end{equation}
We consider two cases: $k$ even and $k$ odd.

If $k$ is even and $k=2t$, from \ref{qee} we deduce that
\begin{equation}
c \equiv \pm p^t \text{ mod } q^2.
\end{equation}
If $p^t \geq \frac{q^2}{2} $, 
\begin{equation}
dist(I,W)=2t \geq \log_p \frac{q^4}{4},
\end{equation}
a contradiction. Consequently, $p^t < \frac{q^2}{2} $. Since $c$ is even, then $c \neq \pm p^t$. We deduce that 
\begin{equation}
c=\pm p^t+lq^2 \text{  for } l\neq 0.
\end{equation}
Therefore,
\begin{equation}
c\geq \frac{1}{2}q^2.
\end{equation}
Hence,
\begin{equation}\begin{split}
p^{2t} \geq \frac{1}{4} q^4,
\\
dist(I,W)=2t \geq \log_p \frac{q^4}{4}.
\end{split}
\end{equation}
This is a contradiction. Therefore $k=2t_0+1$ for some $t_0$.

We now investigate the case where  $$\text{ dist }(I,I^{\prime})=2t+1 < \log_p \frac{q^4}{4}$$  and  $$\text{ dist } (I,W)=2t_0 +1 < \log_p \frac{q^4}{4}.$$ $\text{ dist }(I,I^{\prime})=2t+1 $ gives us a solution to the following diophantine equation
\begin{equation} \label{equa1}
a^2+b^2+c^2+d^2=p^{2t+1} <  \frac{q^4}{4}.
\end{equation} 
Where $b\equiv c \equiv d \equiv 0 \text{ mod } 2q $ and $a \equiv 1 \text{ mod } 2$. At least one of $b, c, d $ is nonzero. Hence
\begin{equation}
4q^2<p^{2t+1}<q^4.
\end{equation}
 $\text { dist }(I,W)=2t_0 +1 < \log_p \frac{q^4}{4}$, gives us a solution to the following diophantine equation
\begin{equation}\label{equa2}
a_0^2+b_0^2+c_0^2+d_0^2=p^{2t_0+1} < q^4/4.
\end{equation} 
Where $b_0\equiv a_0 \equiv d_0 \equiv 0 \text{ mod q}$ and $a_0 \equiv 1 \text{ mod } 2$. From \ref{equa1} and \ref{equa2} we deduce that
\begin{equation}\begin{split} \label{main}
a^2 \equiv p^{2t+1} \text{ mod } q^2 \text{ and $a$ is odd } a<p^{t+1/2} <q^2/2,
\\
c^2  \equiv p^{2t_0+1} \text{ mod } q^2 \text{ and $c$ is even } c<p^{t_0+1/2} <q^2/2.
\end{split}
\end{equation}
Without loss of generality we assume that $t_0 > t$, from \ref{main} we deduce that
\begin{equation}
\pm ap^{t_0 - t}=c.
\end{equation}
However, this is incompatible with the parities of $a$ and $c$. Hence, we conclude Theorem~\ref{ex}.
\end{proof}

\section{ Visiting almost all points after $(1+\epsilon) \log_{k-1}(n) $ steps }\label{visit}
In this section, we show that if we pick two random points from a  $k$-regular Ramanujan graph $G$,  almost surely they have a distance less than 
\begin{equation}
(1+\epsilon)\log_{k-1}(n). 
\end{equation}
The idea is to use the spectral gap of the adjacency matrix of the Ramanujan graphs to  prove an upper bound on the variance. A similar strategy has been implemented by Sarnak ; see \cite[Page~28]{Sarnak3}. 
\begin{proof}
Let $A(x,y)$ be the adjacency matrix of the Ramanujan graph $G$, i.e.
\begin{equation}
A(x,y):=\begin{cases} 1 \text{  if }  x\sim y \\   0 \text  {  otherwise }       \end{cases}.
\end{equation}

Since $A(x,y)$ is a symmetric  matrix, it is diagonalizable. We can write the spectral expansion of this matrix by the set of its eigenfunctions. Namely,
\begin{equation}\label{Selbergg}
A(x,y)=\frac{k}{\|G\|} +\sum_{j} \lambda_j \phi_{j}(x) \phi_j(y),
\end{equation}
where $\big\{ \phi_j \big\}$ is the orthonormal basis of the nontrivial eigenfunctions with eigenvalues $\big\{ \lambda_j \big\}$ for the adjacency matrix $A(x,y)$. Since we assumed that $G$ is a Ramanujan graph, then $|\lambda_j| \leq 2\sqrt{k-1} $. We change the variables and write 

\begin{equation}\label{chang}
\lambda_j=2\sqrt{k-1} \cos \theta_j.
\end{equation}
We define $S(R):=(k-1)^{\frac{R}{2}}U_{R}(\frac{A}{2\sqrt{k-1}})$, where $U_{R}(x)$ is the Chebyshev polynomial of the second kind, i.e.
\begin{equation}
U_{R}(x):=\frac{\sin((R+1) \arccos x)}{\sin(\arccos x)}.
\end{equation}
The following is the spectral expansion of $S(R)$: 
\begin{equation}\label{sss}
S(R)(x,y):= \frac{(k-1)^{\frac{R}{2}}U_{R}(\frac{k}{2\sqrt{k-1}})}{\|G\|} +\sum_{j}(k-1)^{\frac{R}{2}} U_{R}(\frac{\lambda_j}{2\sqrt{k-1}}) \phi_{j}(x) \phi_j(y),
\end{equation}
\begin{rem}\normalfont
Note that if we lift the linear operator $S(R)$ to the universal covering space of the $k$-regular graph $G$, (which is an infinite $k$-regular tree), then $S(R)$ is the linear operator, which takes the average of a function on a sphere with radius $R$. Namely,
\begin{equation}
S(R) f(x) := \sum_{y, \text{dist}(x,y) = R} f(y).
\end{equation}
\\
See \cite[Remark 2]{Sarnak2} for more discussion of this operator.
\end{rem}
\noindent From the formula for the kernel of  $S(R)$ given in \ref{sss}, we obtain
\begin{equation} \begin{split}\label{ident}
S(R )(x,y)= 
\frac{k(k-1)^{R-1}}{\|G\|}+\sum_{j}(k-1)^{\frac{R}{2}} \frac{\sin((R+1) \theta_j)}{\sin\theta_j} \phi_{j}(x) \phi_j(y).
\end{split}
\end{equation}
We calculate the variance over $y$. For $i\neq j$, we have $\sum_{y\in G} \phi_i(y)\phi_j(y)=0$ and $\sum_{y\in G} \phi_i(y)^2=1$, So only the diagonal terms remain in the following summation:
\begin{equation}\begin{split}\label{vv}
\text{ Var }(x):=\sum_{y\in G}\| S(R )(x,y) - \frac{k(k-1)^{R-1}}{\|G\|}  \|^2 
\\
=\sum_{j} (k-1)^{R} \frac{(\sin(R+1) \theta_j)^2}{(\sin\theta_j)^2} \phi_{j}(x)^2.
\end{split}
\end{equation}
Since $\big\{ \phi_j \big\}$ is an orthonormal basis, we have
\begin{equation} \label{trace}
1=\sum_{y\in G} \delta(x,y) dy= \frac{1}{\|G\|}+\sum_{j}  \phi_{j}(x)^2.
\end{equation}
We also have the following trivial trigonometric inequality, which is derived from the geometric series summation formula :
\begin{equation}\label{trig}
|\frac{\sin(R+1) \theta}{\sin \theta}|=|\sum_{j=0}^{R}e^{i\theta}           |\leq R+1.
\end{equation}
From \ref{trace} and \ref{trig}, we obtain
\begin{equation}\label{inq2}
\text{ Var } \leq (R+1)^2 (k-1)^R. 
\end{equation}
We define
\begin{equation}
M:= \left\{ y : S(R)(x,y)=0 \right\}.
\end{equation}
 Note that $M$ is the set of all vertices $y\in G$ , such that there is no path from $x$ to $y$ with length $R$. Therefore, this is exactly the set $M(x,R)$ as defined in the Theorem~\ref{kk}. 
By the definition of the Var given in \ref{vv},
\begin{equation}\label{inq1}
\|M\| |\frac{k(k-1)^{R-1}}{\|G\|}  |^2 \leq \text{ Var}.
\end{equation}
From \ref{inq1} and \ref{inq2},  we have
\begin{equation}
\|M\|\|(k-1)^{R}\| < \|G\|^2(R+1)^2.
\end{equation}
If we choose $R>(1+\epsilon)\log_{k-1}(n)$, then
\begin{equation}
\|M\| \leq n^{1-\epsilon}(1+R)^2.
\end{equation}
Therefore, we  conclude the Theorem~\ref{kk}.

\end{proof}

\section{Numerical Results}\label{data}

In this section, we present our numerical experiments for the diameter of the family of 6-regular LPS Ramanujan graphs $X_{5,q}$ and compare it with the diameter of a family of 6-regular random Cayley graphs $Z^{q}$. Our numerical experiments show that the ratio of the diameter by the logarithm of the number of vertices $\frac{\text{diam}}{\log_5 |X_{5,q}|}$ converges to $4/3$ as $q \to \infty$ for the LPS Ramanujan graphs $X_{5,q}$ .  On the other hand   $\frac{\text{diam}}{\log_5 |Z^{q}|}$ converges to $1$ as $q \to \infty$ for the random Cayley graphs  $Z^{q}$ .
\\

\begin{figure}[t]
\vspace{-0.45cm}
\centering
\raisebox{-0.6cm}{
\includegraphics[width=.54\textwidth]{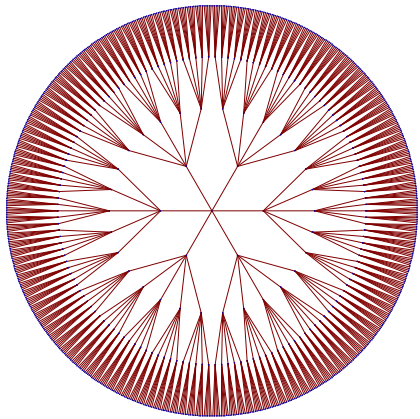}
}\hspace{0.25cm}
\vspace{-0.48cm}
\caption{A ball of radius 4 in the  LPS Ramanujan graphs $X^{5,29}$.}
\label{fig:lps-ball4}
\vspace{-0.2cm}
\end{figure}

We give the detailed construction of the LPS Ramanujan graphs $X^{5,29}$ in what follows. The construction of LPS Ramanujan graphs $X_{5,q}$ requires that $5$ and $-1$ to be quadratic residues mod $q$.  From the reciprocity law we deduce that all the prime factors of $q$ are congruent to 1 or 9 mod 20. The least $q$ with such properties is 29.    We take the integral solutions $\alpha=(x_0,x_1,x_2,x_3)$ of the following diophantine equation 
\begin{equation}\label{gen}
x_0^2+x_1^2+x_2^2+x_3^2=5,
\end{equation}
where $x_0 > 0$  is odd and $x_1,x_2,x_3$ are even numbers. There are exactly $6$ integral solutions with such properties which are listed below:

$$\big\{(1,\pm2,0,0), (1,0,\pm2,0), (1,0,0,\pm2)    \big\}.$$
To each such integral solution $\alpha=(x_0,x_1,x_2,x_3)$, we associates the following matrix in ${\rm PSL}_2(\frac{\mathbb{Z}}{29\mathbb{Z}})$ :
\begin{equation}
\frac{1}{\sqrt{5}} \begin{bmatrix}
x_0+i x_1 & x_2+i x_3
\\
-x_2+ix_3 & x_0 -ix_1
\end{bmatrix},
\end{equation}
where $\sqrt{5}$ and $i$ are the square root of 5 and $-1$ mod 29 respectively. We obtain the following 6 matrices in ${\rm PSL}_2(\frac{\mathbb{Z}}{29\mathbb{Z}})$

$$
S:= \Big\{ \begin{bmatrix} 10 & 0 \\ 0 & 3 \end{bmatrix} , \begin{bmatrix} 3 & 0 \\ 0 & 10 \end{bmatrix},\begin{bmatrix} 8 & 16 \\ 13 & 8 \end{bmatrix},\begin{bmatrix} 21 & 16 \\ 13 & 21 \end{bmatrix},\begin{bmatrix} 21 & 11 \\11 & 21 \end{bmatrix},\begin{bmatrix} 8 &11 \\ 11 & 8 \end{bmatrix} \Big\},
$$
which generate ${\rm PSL}_2(\frac{\mathbb{Z}}{29\mathbb{Z}} )$. The LPS Ramanujan graph $X^{5,29}$ is the Cayley graph of ${\rm PSL}_2(\frac{\mathbb{Z}}{29\mathbb{Z}} )$ with the generator set 
 $S$. The Ramanujan graph $X^{5,29}$ has $12180$ vertices with diameter 8. We note that 
 $$\lceil 4/3\log_5(12180)\rceil=8. $$
%
%
  We show the level structure of  $X^{5,29}$ with root $\begin{bmatrix} 1 & 0 \\ 0 & 1 \end{bmatrix}$ in  table \ref{table}. We note that the girth of this graph is 9
  
 $$\text{girth}(X^{5,29})=9.$$ 
 and this means a ball of radius 4 in the graph $X^{5,29}$ is a tree as illustrated  in Figure \ref{fig:lps-ball4}. For the family of LPS bipartite Ramanujan graphs, Biggs and Boshier determined the asymptotic behavior of the girth of these graphs; see \cite{Biggs}.  They showed that the girth is asymptotic to 
$$(4/3)\log_{k-1}(n).$$

\begin{table}[h!]
\centering
  \begin{tabular}{l|l}
    $r$ & $N(r)$ (Number of vertices of $X^{5,29}$ with distance $r$ from  $\begin{bmatrix} 1 & 0 \\ 0 & 1 \end{bmatrix}$) \\ 
    \hline
    0& 1 \\
    1& 6  \\
    2& 30 \\
    3& 150 \\
    4& 750   \\ 
    5 & 3026  \\
    6& 5970  \\
    7&2195 \\
    8&52\\
  \end{tabular}
  \caption{Level structure of the LPS Ramanujan graphs $X^{5,29}$ \label{table}}
\end{table}

%
%
%

 We give our numerical results for the diameter of the LPS Ramanujan graphs $X_{5,q}$ for  $1 \leq q \leq 229$ in Table \ref{table2}. We note that  $\frac{\text{diam}}{\log_5 n}$ are close to $4/3$. The range for our numerical experiment with the diameter of $X_{5,q}$ is small since the algorithm terminates in $O(q^3)$ operations.  In our very recent work \cite{Naser2}, we developed and implemented a polynomial time algorithm in $\log(q)$ that finds the shortest possible path between diagonal vertices of Ramanujan graphs $X_{p,q}$ under a polynomial time algorithm for factoring and a Cramer type conjecture.  An important feature of our algorithm is that it has been implemented  and it runs and terminates quickly; see \cite[Section 6]{Naser2}.  We give strong numerical evidence that the distance of $W$ from $I$ is asymptotic  to $4/3\log_5(|X_{5,q}|)$ in Table~\ref{alg}. These numerical experiments  are consistent with our conjectures on optimal strong approximation for quadratic forms in 4 variables \cite{Naser}. The conjecture implies that for the LPS Ramanujan graphs $X_{p,q}$ where $p$ is a fixed prime number, the ratio  $\frac{\text{diam}(X_{p,q})}{\log_{p-1} |X_{p,q}| }$ converges to $4/3$ as $q\to \infty$.  Finally, we give our numerical experiments for the diameter of the 6-regular random Cayley graphs ${\rm PSL}_2(\mathbb{Z}/q\mathbb{Z})$. To compare the diameter of the random Cayley graphs with that of the LPS Ramanujan graphs given above, we choose the same set of integers $q$.  We generate 8 random samples  for each $q$, and we give the averaged ratio $\frac{\text{diam}}{\log_5 n}$ in the last column of   Table \ref{table3}.

\begin{table}[h!]
\centering
  \begin{tabular}{  l | l | l | l  }
    
    $q$ & number of vertices of $X_{5,q}$ & Diameter & $\frac{\text{diam}}{\log_5 n}$ \\ 
    \hline
    29& 12180 & 8& 1.36\\
    41& 34440 & 9&1.38 \\
    61& 113460 & 9& 1.24 \\
    89& 352440& 11 & 1.38 \\
    101 & 515100 & 11& 1.34   \\ 
    109 & 647460 & 11 & 1.32  \\
    149& 1653900&  12    &1.34  \\
    181& 3375540& 14 & 1.51\\
    229&6004380 &13&1.34\\
  \end{tabular}
 \caption{LPS Ramanujan graphs $X_{5,q}$ \label{table2}}
\end{table}

\begin{table}[h!]
\centering
  \begin{tabular}{  l | l | l | l  }
    
    $q$ & $d$:= Distance between $W$ and $I$ & $ \frac{d}{\log_5 n}$ \\ 
    \hline
     86028121& 43 & 1.28 \\
    104395301  & 46& 1.35   \\ 
     256203161  & 47 & 1.32  \\
     275604541&  45    &1.26  \\
    472882049 &50&1.36\\
    533000401& 50 & 1.35\\
    613651349  & 50& 1.34\\
    674506081 & 50&1.33 \\
    961748941 & 52& 1.36\\
      32416189381& 57 & 1.28\\
    32416189721& 60 & 1.34\\
    32416189909 & 60 & 1.34\\

  \end{tabular}

  \caption{LPS Ramanujan graphs $X_{5,q}$ \label{alg}}
\end{table}

\

\begin{table}[h!]
\centering
  \begin{tabular}{  l | l | l | l  }
    $q$ & number of vertices of $Z^{q}$ & Diameter & $\frac{\text{diam}}{\log_5 n}$ \\ 
    \hline
    29& 12180 & $8^{\times 6}9^{\times 2}$&1.50\\
    41& 34440 & $9^{\times4}8^{\times 4}$& 1.30 \\
    61& 113460 & $9^{\times 5}10^{\times 3}$& 1.29 \\
    89& 352440& $10^{\times 5}11^{\times3}$ & 1.30 \\
    101 & 515100 & $10^{\times5} 11^{\times3}$& 1.26  \\ 
    109 & 647460 & $10^{\times 4}11^{\times4}$ & 1.26  \\
    149& 1653900&  $11^{\times 6}12^{\times 2}$    &1.25  \\
    181& 3375540& $11^{\times 3}12^{\times 5}$ & 1.24\\
    229&6004380 &$12^{\times 8}$&1.23\\
  \end{tabular}
  \caption{Random Cayley graphs ${\rm PSL}_2(\frac{\mathbb{Z}}{{q\mathbb{Z}}})$ with 6 generators \label{table3}}
\end{table}

($8^{\times 6}9^{\times 2}$ means that 6 of our random samples are 8 and 2 of them are 9). We note that the empirical mean of the ratio $\frac{\text{diam} (Z_q)  }{\log_5 |Z_q|}$ is decreasing in $q$ and  one can easily show    that $$\frac{\text{diam}( Z_q )}{\log_5 |Z_q|}\geq 1.$$ Based on our numerical experiments, we expect that $\frac{\text{diam} (Z_q)  }{\log_5 |Z_q|}$ converges to $1$ in probability as $q \to \infty$ for random Cayley graphs $Z^{q}$ .

\

\

\bibliographystyle{alpha}
\bibliography{revised}

\end{document}